\documentclass[12pt]{amsart}

\usepackage{amscd}
\usepackage{fullpage}
\usepackage{marvosym}
\usepackage{amsmath}
\usepackage{amsgen}
\usepackage{amsbsy}
\usepackage{amsopn}
\usepackage{amsthm}
\usepackage{amscd}
\usepackage{amsfonts}
\usepackage{amssymb}
\usepackage{mathrsfs}
\usepackage{hyperref}
\usepackage{ifsym}
\usepackage{enumerate,enumitem}
 \usepackage{color}

\newcommand{\mc}{\mathcal}

\newcommand{\mb}{\mathbb}

\newcommand{\mf}{\mathfrak}

\DeclareMathOperator{\vol}{Vol}

\newtheorem{Theorem}{Theorem}
\newtheorem{Lemma}[Theorem]{Lemma}
\newtheorem{lemma}[Theorem]{Lemma}
\newtheorem{Corollary}[Theorem]{Corollary}

\newtheorem{Proposition}[Theorem]{Proposition}

\newtheorem{Definition}{Definition}

\newtheorem{Remark}{Remark}

\newcommand{\Vol}{\operatorname{Vol}}
\newcommand{\pa}{\partial}
\newcommand{\del}{\partial}
\newcommand{\op}{\operatorname}
\newcommand{\abs}[1]{\left| #1 \right|}
\newcommand{\norm}[1]{\left\| #1 \right\|}
\newcommand{\inner}[1]{\left\langle #1 \right\rangle}

\newcommand{\set}[1]{\left\{ #1 \right\} }
\newcommand{\tensor}{\otimes}
\newcommand{\grad}{\nabla}
\newcommand{\til}{\widetilde}
\newcommand{\of}{\circ}

\newcommand{\wh}{\widehat}
\DeclareMathOperator{\Aut}{Aut}
\DeclareMathOperator{\Aff}{Aff}
\DeclareMathOperator{\Minvol}{Minvol}
\DeclareMathOperator{\Jac}{Jac}

\newcommand{\R}{{\bf R}}
\newcommand{\N}{{\bf N}}
\newcommand{\Z}{{\bf Z}}

\newcommand{\HH}{{\bf H}}
\newcommand{\HK}{{\bf H}_{\mb{K}}}

\newcommand{\eps}{\epsilon}
\newcommand{\Ga}{\Gamma}

\usepackage{graphicx, amssymb, amsfonts}

\title{On higher graph manifolds}
\author{C.~Connell$^\dagger$ and P.~Su\'arez-Serrato$^\star$}
\thanks{$\dagger$ This work was partially supported by a grant from the Simons Foundation 
\#210442.}
\thanks{$\star$ PSS thankfully acknowledges support for this project from grants IN102716 PAPIIT-DGAPA-UNAM and CN 16-43 CONACyT UC-MEXUS}
   \address{   Department of Mathematics, Indiana University, Bloomington, IN 47405, USA}
              \email{connell@indiana.edu}
    \date{\today}

  \address{Instituto de Matem\'aticas,
              Universidad Nacional Aut\'onoma de M\'exico,
              Ciudad Universitaria, Coyoac\'an,
              04510. M\'exico, D. F.
              M\'exico}

    \email{pablo@im.unam.mx}

\begin{document}

\begin{abstract}
In this short note we {introduce higher graph manifolds} and use a version of the
barycenter technique to characterize when they undergo volume collapse. In the case when
the pure pieces are hyperbolic, we compute the exact value of the minimal volume. We
verify the coarse Baum--Connes conjecture for these manifolds and show that they do not
admit positive scalar curvature metrics. In the case without any pure pieces, we show the
Yamabe invariant vanishes.
\end{abstract}

\maketitle

\section{Introduction}

Consider the class of manifolds we term {\em higher graph manifolds}:
\begin{Definition}
A compact smooth $n$--manifold $M$, $n \geq 3$, is a {\em higher graph manifold}
provided that it can be constructed in the following way:

\begin{enumerate}
\item For every $i = 1, . . . , r$ take a complete finite-volume non-compact
pinched negatively curved $n_i$-manifold $V_i$, where $2 \leq n_i \leq n$. 

\item Denote by $M_i$ the compact smooth manifold with boundary obtained by
``truncating the cusps'' of $V_i$, i.e. by removing from $V_i$ a {(nonmaximal)}
horospherical open neighborhood of each cusp. 

\item Take fiber bundles $Z_i\to M_i$, which are trivial in a neighborhood of $\partial
M_i$ and with fiber an infranilmanifold $N_i$ of dimension $n-n_i$, i.e. $N_i$ is
diffeomorphic to $\widetilde{N_i}/\Gamma_i$, where $\widetilde{N_i}$ is a simply connected
nilpotent Lie group and $\Gamma_i$ is an extension of a lattice $L_i\subset
\widetilde{N_i}$ by a finite group. So the finite cover $N'_i$ of $N_i$ with covering
group $H_i=\Gamma_i / L_i$ is a nilmanifold $\widetilde{N_i} / L_i$. Assume that the
structure group of the bundle $Z_i\to M_i$ reduces to a subgroup of affine transformations
of  $N_i$.

\item Fix a complete pairing of diffeomorphic boundary components between
distinct $Z_i$'s, provided one exists, and glue the paired boundary components
using diffeomorphisms smoothly isotopic\footnote{We rely on some results not
available for smooth pseudo-isotopies} to affine diffeomorphisms of the
boundaries, so as to obtain a connected manifold of dimension $n$.
\end{enumerate}

We will call the $Z_i$'s the {\em pieces} of $M$ and whenever $\dim(M_i)=n$,
then we say $Z_i=M_i$ is a {\em pure piece}.
\end{Definition}

By \cite{KT} the group of connection preserving or affine tranformations,
$\op{Aff}(N)$, of a compact infranilmanifold $N=G/\Gamma$ is precisely 
$$C(G)/(C(G\cap \Gamma)\rtimes \op{Aut}(\Gamma)).$$ 
Here $C(X)$ stands for the center of $X$. It can easily be computed that
$\pi_k(B\op{Aff}(N))$ is always finitely generated, but
$\pi_k(\op{Diff}(G/\Gamma))$ can be infinitely generated in some cases. Hence,
not all bundles above arise with the given structure group condition, apart from
the restriction that they are trivial on the collar.

\begin{Remark}
We point out that the notion of higher graph manifolds defined here includes:
\begin{itemize}[leftmargin=0.25in]
\item The class of {\em generalized graph manifolds} introduced in the recent
paper by Frigerio, Lafont and Sisto \cite{FLS}. In their definition, the pieces
$V_i$ in item $(1)$ above are required to be hyperbolic with toral  boundary
cusps, the $N_i$ in item $(3)$ are required to be tori, and the gluing
diffeomorphisms in item $(4)$ are required to be affine diffeomorphims.

\item The family of manifolds called {\em cusp-decomposable manifolds} studied
by T. Tam Nguyen Phan in \cite{T} where interesting (non)rigid properties are
explored. Seen as higher graph manifolds they have only pure pieces.

\item The {\em affine  twisted doubles} of hyperbolic manifolds. C.S. Aravinda
and T. Farrell study in \cite{AF} the existence of nonpositively curved metrics
for this class of spaces.
\end{itemize}
\end{Remark}

\begin{Remark}
The higher graph manifolds do not strictly generalize the notion of graph
manifold in dimension $3$ since the latter are permitted to have stand-alone
Seifert-fibered pieces such as a $3$--torus. Conversely, a graph manifold in
dimension 3 does not contain hyperbolic pieces that are permitted in these
higher graph 3-manifolds.

In fact F. Waldhausen \cite{W} classified graph manifolds in dimension $3$. The
construction in \cite{W} justified the use of the term {\bf graph}, just like in
our case here. However, we only consider nonsingular bundles.
\end{Remark}

\begin{Remark}
We also point out that in \cite{FLS} the authors produce some simple and
not-so-simple examples of higher graph manifolds that do not admit any CAT(0)
metric. For example, using only pieces that are products of hyperbolic manifolds
and tori they exhibit examples of higher graph manifolds whose fundamental
groups have distorted $\Z^2$ subgroups which CAT(0) groups cannot have.
They also produce more sophisticated examples that fail to be CAT(0) for more
subtle reasons. Indeed, such examples highlight the increased complexity arising
from boundary gluing diffeomorphisms which are isotopic to affine maps instead
of just isometries.
\end{Remark}

Now we present our main results. The minimal volumes $\Minvol$ and $\Vol_{K}$
are closely related smooth topological invariants of $M$ defined as \cite{Gro}, for smooth $C^{\infty}$ metrics $g$,
$$\Minvol(M)=\inf \set{\vol(M,g) \, | \, -1\leq K_g\leq 1}$$
and
$$\Vol_K(M)=\inf \set{\vol(M,g) \, | \, -1\leq K_g},$$
where $K_g$ represents the sectional curvatures of $g$.

\begin{Theorem}\label{thm:zero}
Suppose that $M$ is a higher graph manifold such that there are no pure pieces,
then $\Minvol(M)=\Vol_K(M)=0$.
\end{Theorem}

An important property of $\Minvol$ is that it is capable of distinguishing
between different smooth structures on an underlying topological manifold. This
was first shown by Bessi\`eres in \cite{B} and there are known examples of
manifolds where even vanishing of $\Minvol$ depends on the smooth structure \cite{Kot}.

For a subclass of higher graph manifolds constructed as doubles of hyperbolic
manifolds P. Ontaneda showed the existence of various distinct interesting
smooth structures (some admit nonpositively curved metrics and others do not)
\cite{Ont}.

The simplicial volume $||M||$ of a closed orientable manifold $M$ is defined as
the infimum of $\Sigma_{i}|r_{i}|$ where $r_{i}$ are the coefficients of any
\emph{real} singular cycle representing the fundamental class of $M$. This
invariant was introduced by M. Gromov in order to find lower bounds for the
minimal volume.

For oriented connected compact $n$-manifolds  with boundary, the relative
simplicial volume $\norm{Z_i,\pa Z_i}$ is defined as the $\ell^1$--semi-norm on
relative singular homology $H_n(Z_i,\pa Z_i;{\bf R})$ of the relative
fundamental class $[Z_i,\pa Z_i]\in H_n(Z_i,\pa Z_i;{\bf R})$.

Our second theorem covers the converse to Theorem \ref{thm:zero}.

\begin{Theorem}\label{thm:pos}
For any higher graph manifold $M$, the simplicial volume satisfies
$$\norm{M}= \sum_i \norm{Z_i,\pa Z_i}.$$
In particular it is positive if there is one pure piece.
\end{Theorem}

The {\em volume growth entropy} $h(g)$ of any manifold $(M,g)$ is the quantity
$$h(g)=\limsup_{R\to\infty}\frac{\log \vol B(x,R)}{R},$$
where $B(x,R)$ is the geodesic ball of radius $R$ in {the universal covering}
$\til{M}$. Whenever $M$ covers a compact manifold, the limsup in the
definition of $h(g)$ may be replaced by a limit \cite{M}.

Similarly, we define the {\em minimal volume entropy} of any manifold $M$ as
$$h(M)=\inf \set{h(g)\,|\, \vol(M,g)=1}.$$

Our next theorem gives bounds for $\Minvol$ and $\Vol_K$ {and it relates them to the
minimal volume entropy}.

\begin{Theorem}\label{thm:minent}
Let $M$ be a higher graph manifold all of whose pure pieces $M_1,\cdots,M_k$ are
locally symmetric, and let $M_{sym}=\coprod_{i=1}^k V_i$, where $V_i$ is the
interior of $M_i$. The $V_i$ are all quotients of the same symmetric space, and
if $g_{sym}$ represents the symmetric metric on $M_{sym}$ normalized to have
sectional curvatures between $-4$ and $-1$, then
\begin{equation}\label{Thm3-ineq1}
h_0\vol(M_{sym},g_{sym})^{\frac1n}\leq h(M)\leq 2(n-1)\vol(M_{sym},g_{sym})^{\frac1n}
\end{equation}
where $h_0\in\set{n-1,n,n+2,22}$ depending on whether all of the pure pieces are
of real hyperbolic, complex hyperbolic, quaternionic hyperbolic, or Cayley
hyperbolic type. Moreover,
\begin{equation}\label{Thm3-ineq2}\frac{h_0^n}{(n-1)^n}\vol(M_{sym},g_{sym})\leq \Vol_K(M)\leq \Minvol(M)\leq 
2^n\vol(M_{sym},g_{sym}).
\end{equation}
(If $k=0$ we interpret all quantities as $0$.)
\end{Theorem}

We can give the exact value for $\Minvol(M)$ when the pure pieces are
hyperbolic.
\begin{Theorem}\label{thm:minvol}
With the notation of the previous theorem, when the pure pieces of $M$ are
hyperbolic, the minimal volume of $M$ is precisely
$$\Minvol(M)=\Vol_K(M)=\Vol(M_{sym},g_{sym})=\sum_{i=1}^k \Vol(V_i,g_{hyp})$$
where $g_{hyp}$ represents the complete finite volume constant curvature $-1$
metric on each $V_i$.
\end{Theorem}

This answers a problem posed in \cite{FLS} positively.

\begin{Remark}
Even when the pure pieces are hyperbolic we do not know if $h(M)$ is exactly
achieved by the locally symmetric metric, because the entropy is not necessarily
continuous with respect to the Gromov-Hausdorff topology.
\end{Remark}

Applying results about asymptotic dimension of groups and the fact that higher graph 
manifolds are aspherical, we are able to show:

\begin{Theorem}\label{thm:psc}
Let $M$ be a higher graph manifold. Then $\pi_1(M)$ has finite asymptotic 
dimension and $M$ does not admit any smooth metric of positive scalar curvature.
\end{Theorem}

The finite asymptotic dimension of $\pi_1(M)$ has several consequences:

For a group $G$ equipped with a left-invariant metric, $d$, consider the reduced
$C^{\ast}$--algebra $C_{r}^{\ast}(G)$. We say that $G$ satisfies the coarse
Baum--Connes conjecture if the coarse assembly map (a homomorphism) of locally
finite complex $K$--homology of the the proper metric space $(G,d)$ to the
algebraic $K$--theory of $C_{r}^{\ast}(G)$ is an isomorphism.

It follows from the work of G. Yu in \cite{Y} that the coarse Baum--Connes
conjecture holds for higher graph manifolds.

\begin{Corollary}\label{cor:bc}
Let $M$ be a higher graph manifold, then {$\pi_1(M)$} satisfies the coarse 
Baum--Connes conjecture.
\end{Corollary}

We take this opportunity to record two corollaries that follow from the validity
of the coarse Baum--Connes conjecture for these groups (see \cite{Y} and
\cite{V}). For a connected, closed, oriented, smooth manifold $M$ its Hirzebruch
$L$--class $L_{M}$ is a polynomial  in the Pontrjagin classes of $M$. The
Novikov conjecture states that, for the map $f:M\to B(\pi_1M)$ that classifies
the universal cover, the class $f_{\ast}(L_{M})\in H_{\ast}(B(\pi_1M), {\bf Q})$
is an oriented homotopy invariant. Equivalently, the higher signatures
$\sigma_{u}(M)\in {\bf Q}$, defined by $\sigma_{u}(M)=\langle L_{M} \cup
f^{\ast}(u), [M]\rangle$ are, for all $u\in H^{\ast}(B(\pi_1M), {\bf Q})$,
oriented homotopy invariants.

\begin{Corollary}\label{cor:nc}
Let $M$ be an oriented higher graph manifold, then the Novikov
conjecture holds for $M$.
\end{Corollary}

In fact, by Yu \cite{Y}, Corollary  \ref{cor:bc} actually implies that such a
manifold has a {\em zero in the spectrum:}

\begin{Corollary}\label{cor:zerospec}
Let $M$ be an oriented compact higher graph manifold, then there exists a $p\geq
0$, such that zero belongs to the spectrum of the Laplace-Beltrami operator
$\Delta_{p}$ acting on square-integrable $p$--forms on the universal cover of
$M$.
\end{Corollary}

Consider a fixed conformal class $\gamma $ of smooth metrics $g$ on the smooth
closed manifold $M$, and let the Yamabe constant of $(M, \gamma )$ be
$$ \mathcal{Y}(M, \gamma )= \inf\limits_{g \in \gamma} \frac{\int\limits_{M} s_{g} 
dvol_{g}}{(\Vol(M,g))^{1-2/n }}.$$
Here $s_{g}$ represents the scalar curvature of $g$. The \emph{Yamabe invariant}
is then defined to be $ \mathcal{Y}(M) = \sup \mathcal{Y}(M, \gamma ),$ where
the supremum is taken over all conformal classes of metrics on $M$. As a
consequence we derive that the Yamabe invariant vanishes when all the pieces
have infranilpotent fibers of positive dimension:

\begin{Corollary}\label{cor:yam}
Let $M$ be a compact higher graph manifold without pure pieces. Then the Yamabe
invariant of $M$ equals 0.
\end{Corollary}

The proofs of all of the results are found in the following section.

\medskip {\bf Acknowledgements:} PSS thanks support from PAPIIT, Universidad
Nacional Aut\'onoma de M\'exico, from the CNRS in France and also LAISLA (CNRS-CONACyT joint project). We benefited from
UNAM's international academic exchange programme for a visit to Mexico City by
the first author and to Bloomington by the second author. We warmly thank the
Laboratoire Jean Leray in Nantes for their gracious hospitality during the final
writing stages, and in particular Gilles Carron for comments on a previous
version. 

\section{Proofs}

In our following proof of Theorem \ref{thm:zero} we will adapt arguments first
explained by K. Fukaya in \cite{F}, where he constructed an explicit sequence of
volume collapsing metrics on fibre bundles with infranilpotent fibres and affine
structure group.

Briefly, Fukaya's insight was to use the Lie algebra structure to organize
charts of $Z_{i}$ and describe a collapsing sequence of metrics using the level structure
of the Lie algebra $\mathfrak{n}$ of the nilpotent Lie group $\widetilde{N}$
that covers the fiber $N$. First we choose a connection relative to the
structure group of the fibration. A metric is then defined using the horizontal
and vertical subspaces set by the connection and, crucially, the central
direction of $\mathfrak{n}$, which is the part of the metric that collapses. As
the structure group is affine, it is possible to construct charts compatible
with this metric, which moreover keep the collapsing directions of the metric
consistent (since they come from central elements, and the charts are glued
using affine transformations). So we obtain a parameterized family of smooth
metrics on all of $Z_{i}$ whose volume vanishes in the limit (in fact the base
$M_{i}$ is the Gromov-Hausdorff limit). The curvature is shown to be bounded
throughout this process.

For the proof we will first need the following definition and technical lemma which may be of 
interest in its own right.

\begin{Definition}\label{def:rel_c_bdd}
We say that two $C^2$ metrics $\rho$ and $\sigma$ are {\em $C$-relatively bounded up to first
order} if there is an atlas for which $\norm{\rho_{ik}\sigma^{kj} -\delta_{ij}}_{C^1}<C$
and $\norm{\sigma_{ik}\rho^{kj} -\delta_{ij}}_{C^1}<C$  for all $i,j$ and for some
constant $C<\infty$. (Note that this is equivalent to  $\norm{\rho}_{C^1,\sigma}<C$ and
$\norm{\sigma}_{C^1,\rho}<C$ for a possibly different $C$).
\end{Definition}

\begin{lemma}\label{lem:curv_of_sums}
	Suppose that the $C^2$ metrics $\rho$ and $\sigma$ are $C$-relatively bounded up to first
	order, then, for $\alpha$ in $[0,1]$, the sectional curvatures of $g=\alpha\sigma+(1-\alpha)\rho$ are bounded above and below in
	terms of those of $\rho$ and $\sigma$ and independently of $\alpha$.
\end{lemma}

\begin{proof}
Choose normal coordinates for the metric $g=\alpha \sigma +(1-\alpha)\rho$ at a fixed point 
$p$ so that $g_{ij}=\delta_{ij}=\alpha \rho_{ij} +(1-\alpha)\sigma_{ij}$, and the assumptions 
imply that $\rho_{ij}$ and $\sigma_{ij}$ are both uniformly bounded above and away from 0 at 
$p$. The norm conditions imply that $\rho_{ij,k}$ and $\sigma_{ij,k}$ are bounded above in 
these coordinates. 

It follows that the Christoffel symbols 
$^{\rho}\Ga^{i}_{jk}=\rho^{is}(\rho_{sj,k}+\rho_{sk,j}-\rho_{jk,s})$ are bounded and 
likewise for $^{\sigma}\Ga^{i}_{jk}$. 

The formula for the curvature tensor for $\rho$ is:
\[
{}^{\rho}R^{l}_{ijk}= {}^{\rho}\Ga^{l}_{ik,j}- {}^{\rho}\Ga_{ij,k}^l+ 
{}^{\rho}\Ga_{ik}^s{}^{\rho}\Ga_{sj}^l- {}^{\rho}\Ga_{ij}^s {}^{\rho}\Ga_{sk}^l.
\]

We consider each terms of $^{\rho}\Ga^{l}_{jk,i}$ separately. It suffices to consider 
$\rho^{lm}\rho_{mj,k}$. We have
\[
\del_i \rho^{lm}\rho_{mj,k}= \rho^{lm}\rho_{mj,k,i} - \rho^{lr}\rho_{rs,i}\rho^{sm}\rho_{mj,k}
\]
Since the latter terms are bounded, and the Christoffel symbols are bounded and $\rho^{lm}$ is 
bounded, we may write:
\[
(\rho_{mi,k,j}+\rho_{mk,i,j}-\rho_{ik,m,j})-(\rho_{mi,j,k}+\rho_{mj,i,k}-\rho_{ij,m,k})=\rho_{ml}
 {}^{\rho}R^{l}_{ijk}+{}^{\rho}C_{mijk}
\] 
where the term ${}^{\rho}C_{mijk}$ involves only sums of terms of the form
$\rho_{ms,i}\rho^{sl}\rho_{lj,k}$ and is therefore uniformly bounded. The same is true for
the analogous $\sigma$ terms.

The normal coordinates also imply, $g_{ij,k}=0$ at $p$. So the curvature tensor of the metric 
$g_{ij}=\alpha \rho_{ij} +(1-\alpha)\sigma_{ij}$ is 
\[
R^{l}_{ijk}=(g_{li,k,j}+g_{lk,i,j}-g_{ik,l,j})-(g_{li,j,k}+g_{lj,i,k}-g_{ij,l,k})
\]
Since $g_{li,k,j}=\alpha \rho_{li,k,j}+(1-\alpha)\sigma_{li,k,j}$ and similarly for the other 
terms, we may write this as the sum 
\[
R^{l}_{ijk}=\alpha \rho_{lm} {}^{\rho}R^{m}_{ijk}+(1-\alpha)\sigma_{lm} 
{}^{\sigma}R^{m}_{ijk}+ \alpha {}^{\rho}C_{mijk}+ (1-\alpha){}^{\sigma}C_{mijk}
\]
where $\abs{\alpha {}^{\rho}C_{mijk}+ (1-\alpha){}^{\sigma}C_{mijk}}$, $\abs{\rho_{lm}}$ 
and $\abs{\sigma_{lm}}$ are all bounded above uniformly and independently of $\alpha$ as 
well.

Finally we note that the choice of the point $p$ was arbitrary.
\end{proof}

Now we are ready for the proof of the first Theorem.

\begin{proof}[Proof of Theorem \ref{thm:zero}]

The explicit construction of a volume collapsed sequence of metrics for fibre
bundles with infranilpotent fibre of Theorem 0-7 in \cite{F} applies directly to
each of the pieces $Z_{i}$ that make up $M$. It is enough to show that these
metrics can be used to provide a sequence of volume collapsing metrics on the
whole of $M$.

For an infra-nilmanifold $N$, we may write $N=\til{N}/\Ga$ where $\til{N}$ is a
nilpotent Lie group and $\Ga<\Aff(\til{N})=\til{N}\rtimes \Aut(\til{N})$ with
$[\Ga\cap \til{N}:\Ga]$ finite.

We now consider two pieces reindexed as $Z_1$ and $Z_2$ that will be joined across a 
boundary component that we also index with $i=1,2$ respectively. 

Fukaya's construction gives a volume collapsing sequence of metrics
$g_{i,\epsilon}$ of bounded curvature on $Z_{i}$ for which the infranil fiber
$N_{i}$ has diameter $\eps$ and which are left invariant under the covering
nilpotent group, denoted $\til{N}_{i}$. Moreover, the pointed Gromov-Hasudorff
limit metric $g_{i,0}=\lim_{\eps\to 0}g_{i,\eps}$ will be any choice of given metric on
the base $M_i$.

For the metric on $M_i$, first consider horospherical parametrizations of the truncated 
portion of the cusps of $V_i$ extending from $\partial M_{i}$, used in the construction of
$Z_i$, so that the negatively curved metrics in these parameters are described
by a generalized warped Riemannian product on a neighborhood of each boundary
component:
$$dy^2_{c_{i}(t)}+dt^2$$
Here $dy^2_{c_{i}(t)}$ represents the metric on the cross-section at parameter $t$,
$c_{i}(t) \cong C_{i}\times\set{t}$, of the cusp of $V_{i}$ topologically identified with
$C_{i}\times [0,\infty)$ where the portion of the cusp of $V_i$ that was truncated to
obtain $M_i$ corresponds to $t\in (0, \infty)$. 

We now change this metric on $C_{i}\times [t,\infty)$ to a new metric denoted by
$dy_{i,t}^2+dt^2$ where $dy_{i,t}$ is defined as follows. Let $\mf{c}_i$ denote
the Lie algebra of $\til{C}_i$. This has an associated flag decomposition from
the nilpotent structure $\mf{c}_i=L_1\supset L_2 \supset \dots \supset
L_{k+1}=0$ such that $L_{i+1}=[L_1,L_i]$ for all $i={1,\dots,k}$. Choose any
compact subgroup $K<\Aut(\til{C}_i)$ containing $\Ga\cap \Aut(\til{C}_i)$. Now
we choose any $K$-invariant metric $\inner{,}_0$ on $\mf{c}_i$ and define
\[
F_i=\set{X\in L_{i}: \inner{X,Y}_0=0 \text{ for all } Y\in L_{i+1}}.
\]
It follows that $\mf{c}=F_1\oplus F_2 \oplus \cdots \oplus F_k$. Let $s$ be a fixed constant. Now we define
for all $t\geq 1$ a $K$-invariant inner product $\inner{,}_t$ on $\mf{c}_i$ by
$\inner{X,Y}_t=s^2e^{-2it}\inner{X,Y}_0$ for $X,Y\in F_i$ and  $\inner{X,Y}_t=0$
for $X\in F_i$ and $Y\in F_j$ when $i\neq j$. In particular for any constant
$s>0$ this defines a $\til{C}_i\rtimes K$-invariant metric on $\til{C}_i$ which
descends to the metric $dy^2_{i,t}$ on $C_i$. Theorem 11.1 of \cite{BK2} shows
that for all sufficiently large constants $s>0$ the $\til{C}_i\rtimes
K$-invariant metric $dy^2_{i,t}+dt^2$ on $C_i\times [1,\infty)$ has pinched
curvatures lying in $[-(k+1)^2,-1]$.

We smoothly extend this metric to $C_i\times [0,\infty)$ so that
$dy^2_{i,t}=dy^2_{c_i(t)}$ for $t\in [0,\frac12]$ and the curvatures remain
uniformly bounded. We can do this since the metrics $dy^2_{c_i(t)}+dt^2$ and
$dy^2_{c_i(t)}+dt^2$ have uniformly bounded curvatures and injectivity radius in
collar neighborhoods corresponding to $t\in[0,\frac12]$ and $t\in [1,2]$
respectively. For example, we could explicitly define
$dy^2_{i,t}=(1-\beta(t))dy^2_{c_i(\frac12)}+\beta(t)dy^2_{i,1}$ for
$t\in[\frac12,1]$ for any smooth function $\beta$ supported in $[\frac12,1]$ and
increasing from $0$ to $1$ with all derivatives vanishing at $t=\frac12$ and
$t=1$.

Since the $N_i$ bundle on the cusp was assumed to be a product $N_i\times 
C_i\times[0,\infty)$, the structure group restricted to the cusp is trivial. Hence in 
Section 6 of \cite{F} we may choose the bundle connection so that it is parallel on the 
cusp and therefore the resulting metric $g_{i,\eps}$ takes the form of a product metric 
\[
g_{i,\eps}=dx^2_{i,\eps} + dy^2_{i,t}+dt^2
\]
for some fixed left invariant metric $dx^2_{i,\eps}$ of diameter $\eps$ on $N_i$.

\begin{lemma}\label{lemma:bilip}
%If $t_i=t_i(\eps)$ is chosen so that $N_i$ and $C_i$  have the
%same diameters with respect to $dx^2_{i,\eps}$ and $dy^2_{i,t}$ respectively and for $i=1,2$, 
Setting $t_{i} = -\log(\eps) + \log(s)$ then $f:N_1\times C_1 \to N_2\times C_2$ is $C$-bilipschitz for some 
constant $C>1$ independent of $\eps$.
\end{lemma}

\begin{proof}
%If $\eps$ is sufficiently large so that $t_i<1$ then there is nothing to prove, so we assume $t_i>1$. 

The metrics above on $C_i$ and the metrics on $N_i$ constructed in the proof of
Theorem 0.7 of \cite{F} are defined identically once the parameters are matched
as $\eps=se^{-t}$. Consequently, for $t_1=t_2=-\log \eps+\log s$ the combined metric
$dx^2_{1,\eps}+ dy^2_{1,t}$ on $N_1\times C_1$ has size $\eps^{i}$ on the left
translate of the subspace $F_i$ relative to the choice of initial metric
$\inner{,}_0$ on the combined algebra $L_1=\mf{n}_1\oplus\mf{c}_1$. The same
holds for $dx^2_{2,\eps}+ dy^2_{2,t}$ relative to the choice of initial metric 
$\inner{,}_0'$ on $L_1'=\mf{n}_2\oplus\mf{c}_2\cong L_1$. In particular, up to an 
additive constant depending only on $\inner{,}_0$, this choice of $t_i$ will make the 
diameters of $C_i$ and $N_i$ both $O(\eps)$.

Since $N_1\times C_1$ is isomorphic to $N_2\times C_2$ and the diffeomorphism
$f$ is affine, its derivative sends the levels $L_i$ into $L_i'$, even though it
may mix factors from between $C_i$ and $N_i$. While $df$ may not send the
subspace $F_i$ precisely onto $F_i'$, as they depend on the choice of initial
metrics $\inner{,}_0$ and $\inner{,}_0'$ on their respective sides, nevertheless
$df$ does send $F_i$ into $\oplus_{j\geq i}F_j'$. Moreover, since $df$ cannot
take $L_i$ into $L_{i+1}'$ the linear matrix coefficients of $df$ corresponding
to the $F_i'$ subspace of the image only depend on $\inner{,}_0$ and
$\inner{,}_0'$ and so are uniformly bounded and bounded away from 0 independent
of $\eps$.

Hence, for the choice of $t_i(\eps)$ above, $df$ is uniformly bounded and
bounded away from 0 on each $F_i$, and hence on all of $L_1$. Therefore with
respect to these two metrics, $f$ is uniformly Lipschitz independent of $\eps$.
\end{proof}

Now for $i=1,2$ let $t_i=t_i(\eps)$ be the time parameter from the lemma such that 
$c_i(t_i)$ has diameter $\eps$.

Fix a constant $s$ and let $A^i_s$ be the portion of the cusp of $V_i$ corresponding to $t\in [0,s]$. Now without loss of generality we may assume $t_i>1$.  Observe that we may smooth the warped product metric on $A^i_{t+R}$, for $t\in [t_i-1,t_i]$ and $R>0$,  to the fixed 
product metric 
$$dy_{c_i(t_i)}^2+ dt^2$$
 for all $t\in [t_i,t_i+R]$.

Let $\beta:\R\to[0,1]$ be a smooth bump function with values:
\begin{eqnarray*}
\beta(t)=1 & \quad {\rm for} \quad & t \in (-\infty,t_1]\\
\beta(t)=0 &  \quad {\rm for}  \quad& t \in [t_1+1,\infty)\\
\end{eqnarray*}

Now create the adjunction space 
\[
A=N_1\times A^1_{t_1+1}\cup_q N_2\times A^2_{t_2+1}.
\]
Here $q$ is the diffeomorphism on the product $N_1\times (A^1_{t_1+1}\setminus A^1_{t_1})$
and $N_2\times (A^2_{t_2+1}\setminus A^2_{t_2})$ whereby $(x,y,t)\in N_1\times
A^1_{t_1+1}$ is glued to $(f(x,y),t_2+t_1+1-t)\in N_2\times A^2_{t_2+1}$ for all $t\in
[t_1,t_1+1]$.

Note that $Z_1\cup_f Z_2$ is diffeomorphic to $Z_1\cup A\cup Z_2$. On $A$ we put
the warped product metric $g_{A,\eps}$ which on the fiber over
$t\in[t_1,t_1+1]$, following the $A^1$ parametrization, is:
\[
g_{A,\eps}(t)=\beta(t)g_{1,\eps}(t) + (1-\beta(t))f^*g_{2,\eps}(t_2+t_1+1-t),
\]
and is otherwise $g_{1,\eps}$ on $A^1_{t_1}\subset A$ and $g_{2,\eps}$ on 
$A^2_{t_2}\subset A$.

We now show the boundedness of the curvatures of the metrics $g_{A,\eps}$ on $A$, which 
isometrically attach to the metrics $g_{i,\eps}$ at the boundary components of $A$. 

Since $f$ is $C$-biLipschitz, we have that $f^*g_{2,\eps}$ and $g_{1,\eps}$ are
$C$-relatively bounded up to first order for some constant $C\geq 1$. From Lemma
\ref{lem:curv_of_sums} we obtain that the fiberwise metrics $g_{A,\eps}(t)$ defined above
have bounded curvartures independent of the values of $\beta(t)$ since the metrics
$g_{i,\eps}$ have bounded curvature, and $f^*g_{2,\eps}$ is isometric to $g_{2,\eps}$.
Hence the curvatures of $g_{A,\eps}(t)$ are bounded independent of $\eps$. It remains to
bound the curvatures of $g_{A,\eps}$ that involve the $t$ direction.

Suppressing the fixed $\eps>0$, let $g_A(t)$ denote the combined metric on the
fiber over $t$, which we recall is topologically $N_{i}\times C_{i}$, and let
$g$ be the metric on $Z_1\cup A \cup Z_2$. Now, denote the field of the $t$ coordinate as
$\del_t$, and choose any $g_A(t)$-orthogonal geodesic frame $\set{X_i}$ in a
neighborhood of a point $p$ in each fiber over $t$. Writing
$h_i(t)=\sqrt{g_A(t)(X_i,X_i)}$, the fields $Y_i=X_i/h_i$ are $g_A(t)$-orthonormal.
Then we can express the curvature tensor in these coordinates by (see Appendix C
of \cite{Be} correcting \cite{BK}):
\begin{align}\label{eq:warped}
\begin{split}
R_g(Y_i,Y_j,Y_j,Y_i)&= R_{g_A(t)}(Y_i,Y_j,Y_j,Y_i)-\frac{h_i'h_j'}{h_ih_j}\quad \text{if }i\neq j\\
R_g(Y_i,Y_j,Y_k,Y_l)&= R_{g_A(t)}(Y_i,Y_j,Y_k,Y_l) \quad \text{if }\set{i,j}\neq 
\set{k,l}\\
R_g(Y_i,\del_t,\del_t,Y_i)&=-\frac{h_i''}{h_i} \quad \text{and} \quad 
R_{g}(Y_i,\del_t,\del_t,Y_j)=0 \quad \text{if }i\neq j  \\
2R_g(\del_t,Y_i,Y_j,Y_k)&=\inner{[Y_i,Y_j],Y_k}\left(\log \frac{h_k}{h_j}\right)' 
+\inner{[Y_k,Y_i],Y_j}\left(\log \frac{h_j}{h_k}\right)'\\
&\quad  +\inner{[Y_k,Y_j],Y_i}\left(\log \frac{h_i^2}{h_j h_k}\right)'
\end{split}
\end{align}
 
If we let $\set{X_i}$ be a fiberwise orthogonal base for $g$ then each $X_i$ can
be expressed linearly in terms of an orthogonal basis for $g_{1}(t)$, $g_{2}(t)$ and
$\beta(t)$. If we now specialize to a specific bump function like 
$\beta(t)=1-e^{-\frac{1}{1-x^2}}$, then we note that 
\[
\frac{h_i'(t)}{h_i(t)}=\frac12 (\log h_i^2)'(t)=\frac{1}{2 
h_i^2(t)}\left(\beta'(t)(a_1(t)+a_2(t))+\beta(t)(a_1'(t)-a_2'(t))+a_2'(t) \right),
\]
where $a_1(t)=g_{1}(X_i,X_i)$ and $a_2(t)=g_{2}(X_i,X_i)$. Note that $a_1(t)$,
$a_2(t)$ and all their derivatives are bounded in $t$, independently of $\eps$.
Moreover,  $\beta$ and its derivatives are all bounded by $1$, and $h_i^2(t)\geq
\min\set{a_1(t),a_2(t)}$ and that $a_1(t)$ and $a_2(t)$ are $c$-comparable for
some $c>1$. Hence the numerator above is at most $c$ times the denominator, and
so each $\frac{h_i(t)}{h_i(t)}$ is bounded. Taking one more derivative to obtain
$(h_i^2)''$, the same argument shows that $\frac{(h_i^2)''}{h_i^2}$ is bounded
and hence
\[
\frac{h_i''}{h_i}=\frac{(h_i^2)''}{2h_i^2}-\left(\frac{h_i'}{h_i}\right)^2
\] 
is bounded as well. Hence all curvatures are bounded in terms of $R_{g_A(t)}$
independently of $\eps$. Therefore all curvatures are  bounded independent of
$\eps$.

In the end, the resulting metric may only be $C^2$ since the cuspidal Busemann
functions in finite volume manifolds, which we used for the parameter $t$, are
only $C^2$ in general. However, by standard approximation theory (see
\cite{Hi}), any $C^2$ atlas is compatible with a $C^\infty$ one. In this atlas
we can approximate the given metric $g_\eps$ with a $C^\infty$ metric
$g_{\eps,\eta}$ with $C^2$ norm $\norm{g_\eps-g_{\eps,\eta}}_{C^2}$ sufficiently
small so that the difference of both the sectional curvatures and the volumes
are less than $\eta$. Then choosing the bound $\eta=\eps$, then the sequence of
metrics $g_{\eps,\eps}$ is $C^\infty$ and also collapses.
\end{proof}

\begin{Remark}
We note that in the construction above, if the gluing map mixes a portion of the
$C_1$ factor into $N_2$ and $N_1$ into $C_2$, then the diameter of the $A$ portion
tends to infinity throughout the collapse. More specifically, the $A_{t_i}^i$ pieces 
will limit as $t_i\to\infty$ to a cusp neighborhood of the manifold $V_i$.

\end{Remark}

\begin{Remark}
One consequence of the above theorem together with the work of Cheeger, Fukaya
and Gromov is that there exists a global $\mc{F}$--structure (and an
$\mc{N}$--structure) on higher graph manifolds without pure pieces. Moreover, it
seems very plausible that an $O(n)$-invariant lifting of the local infranil
fibrations to the frame bundle would provide one such (mixed) $\mc{N}$-structure, and
that it would even be polarized since the action on $M$ is everywhere locally
free. By the work of Cai and Rong \cite{CR} a polarized $\mc{N}$-structure
implies volume collapse. However, there are a number of subtleties involved in proving 
that the obvious candidate structure actually is one, and so we chose the more 
constructive approach above. 
\end{Remark}

\begin{Remark} We thank an anonoymous referee for pointing out a result of Fukaya in 
\cite{F0}, relevant to the previous proof. The main result of  \cite{F0} implies that a 
sequence of volume collapsing metrics such as the one constructed in the proof above can 
have bounded diameters only if the manifold is diffeomorphic to an infranilmanifold 
fibration. Hence the need to allow the diameters of the metrics constructed in the proof 
to tend to infinity, therefore yielding a volume collapsing sequence of smooth metrics. 
\end{Remark}

\begin{proof}[Proof of Theorem \ref{thm:pos}]
For any boundary component $E$ of any piece $Z_i$, we claim that the homomorphism
$\pi_1(E)\to \pi_1(Z_i)$ induced by inclusion is injective. To see this, we first note
that since $V_i$ has finite volume, the inclusion of the corresponding boundary
component, $B\subset \partial M_i$, is $\pi_1$ injective in $\pi_1(M_i)$. Second, since
the pair $(Z_i,E)$ is the total space of an $N_i$ fibration over $(M_i,B)$, the
injectivity can be deduced from the long exact sequence.

Moreover, the inclusion is amenable simply because $\pi_1(E)$ is amenable. Gromov's {\it
cutting-off} theorem of \cite{Gro} has been extended by T. Kuessner in his thesis
\cite{K} (see also \cite{K2}) to exactly cover the case in question. By Lemma 11 (i),(ii) of \cite{K}, we
obtain that $\norm{M,\pa M}\geq \norm{Z_i,\pa Z_i}+\norm{ Q,\pa Q}$ where $Z_i$ is the
one piece and $Q$ is the complementary manifold. (Lemma 11 (ii) covers the case where two
boundary components of $M_i$ are glued together.) Inductively, we obtain $\norm{M}\geq
\sum_i \norm{Z_i,\pa Z_i}$.

For the reverse inequality, we may apply Lemma 12 in the same thesis since the
fundamental group of the entire collar is amenable.

In \cite{Gro} Gromov proved the positivity of $\norm{M_i,\pa M_i}$ whenever
$\op{int}(M_i)=V_i$ admits a finite volume metric of pinched negative curvature. The non
pure pieces collapse and therefore the simplicial volume is at least the sum of the
relative simplicial volumes of the pure pieces which is positive.
\end{proof}

The following proposition generalizes part of Theorem 4 of J. Souto \cite{S}. We will 
follow the same general approach though significant changes are required at some points. 
We will need this preliminary result to prove Theorems \ref{thm:minent} and 
\ref{thm:minvol}.

\begin{Proposition}\label{prop:entropy}
Let $M$ be a closed higher graph $n$-manifold {whose pure pieces are all locally symmetric
of the same type.} If $(M_{sym},g_{sym})$ represents a locally symmetric metric on the
pure pieces scaled to have curvatures between $-1$ and $-4$, then
$$h(M)\geq h_0\left(\vol(M_{sym},g_{sym})\right)^{\frac1n}$$
where $h_0\in\set{n-1,n,n+2,22}$ depending on whether all of the pure pieces are of real
hyperbolic, complex hyperbolic, quaternionic hyperbolic, or Cayley hyperbolic type.
\end{Proposition}

\begin{proof}
We first consider a CW-complex $X$ defined as the quotient space formed by first
taking $M$ and collapsing each embedded non-pure piece $Z_i\subset M$ to a
point, say $v_i$, with $v_i=v_j$ whenever $Z_i$ and $Z_j$ are connected through
a sequence of non-pure pieces. Additionally, if for each pair of adjacent pure
pieces $M_i$ and $M_j$ we denote the connected components of $\pa M_i\cap \pa
M_j\subset M$ by $\set{B_{ij}^1,\dots, B_{ij}^m}$, we collapse each $B_{ij}^k$
to a distinct point $e_{ij}^k$. The space $X$ has a naturally associated graph
structure $\mc{G}$ whose vertices are the collection of pure pieces together
with the points $\set{v_i}$. The edges of $\mc{G}$ consist of all pairs of
vertices where one vertex corresponds to a pure piece and the other to a vertex
in $\set{v_i}$ lying on the boundary of the pure piece, together with all the
elements of $\set{e_{ij}^k}$ correspond to the edges between pure vertices.

Let $f:M\to X$ be the quotient map. Note that
$$f_*[M]=[M_1,\pa M_1]+\cdots+[M_k,\pa M_k]\in H_n(X,\Z).$$
Let $\mc{S}=\set{v_i}\cup\set{e_{ij}^k}$ be the singular set of collapsed points. We 
observe that the complementary regular set $\mc{R}=X-\mc{S}$ is homeomorphic to 
$\coprod_{i=1}^k V_i$ where the $V_i$ are the interiors of the pure pieces $M_i$.

We note that the usual quotient metric is not locally CAT(-1) since for all
sufficiently small $\delta>0$, among all pairs of points in the same piece $M_i$
at distance $\delta$ from a collapsed point $p\in\mc{S}$ we can always find a
pair $x,y$ whose distance in $M_i$ is $2\delta$. Thus, there are at least two
distinct geodesics joining $x$ and $y$, one traveling through $p$. Instead, we
change the metric in a neighborhood of the points of $\mc{S}$ as we now describe.

Starting with locally symmetric metrics on each $V_i$, we can smoothly deform the metric
in a horospherical neighborhood $U$ of each end with an incomplete warped product metric
which is completed by a single cusp point. To construct this, we note that the initial
metric on $U$ has the warped product form $ds^2=d(\delta(t)g_o)^2+dt^2$ where $g_o$ is the
standard left-invariant Riemannian inner product on the nilpotent Lie algebra ${\mf n}$ of
Heisenberg type which gives rise to the induced metric on the nilmanifold corresponding to
the horospherical cross-section of the end, and $\delta(t)$ is the homothety of ${\mf
n}=V_1\oplus V_2$ which dilates by $t$ on the $V_1$ distribution and by $t^2$ on the
central $V_2=[V_1,V_1]$ subspace. Now we note that the decreasing convex function
$e^{3-\frac{3}{1-t}}$ for $t\in [2/3,1]$ can be smoothly interpolated with $e^{-t}$ for
$t\in[0,1/3]$ to obtain a positive decreasing convex function $f:[0,1]\to [0,1]$.  The
warped product metric we want is $ds^2=d(\delta(f(t))g_o)^2+dt^2$, which from equations
\eqref{eq:warped} we see still has sectional curvatures in $(-\infty,-1]$.

We complete the metric to obtain the desired metric on all of $X$ by adding in
the single points at the finite cusp ends, namely the points of $\mc{S}$. This
completed metric is again not locally CAT($-1$) at points of $\mc{S}$ since there
are arbitrarily short (contractible) geodesic bigons around the cones ending at each $p\in
\mc{S}$. Though $\mc{R}$ is smooth and negatively curved, its universal cover $\til{\mc{R}}$ 
is not geodesic since geodesics between distinct fundamental domains of the cusps want to pass through the missing points corresponding to $\mc{S}$.

Now define $\wh{X}$ to be the quotient space of $\til{M}$ formed by collapsing
each connected component of the lifts of $f^{-1}(\mc{S})$ in $\til{M}$ to
points, and denote the quotient map by $\wh{f}$. Denote these collapsed points
by the countable discrete set $\wh{\mc{S}}$ and the complement of $\wh{\mc{S}}$
in $\wh{X}$ by $\wh{\mc{R}}$. Each component of $\wh{\mc{R}}$ will be
homeomorphic to some component of the universal cover of the regular set
$\til{\mc{R}}$. We metrize each component of $\wh{\mc{R}}$ to be isometric to
its corresponding component of the negatively curved manifold $\til{\mc{R}}$.
The metric completion of each component can be formed by adding on a subset of
the points of $\wh{\mc{S}}$. This endows a metric on $\wh{X}$ by having the
distance between two points in neighboring components of $\wh{\mc{R}}$ to be the
sum of the (finite) distances to and between the unique points in $\wh{\mc{S}}$
that connect the two components in $\wh{X}$.

We claim that $\wh{X}$ a CAT($-1$) space. Observe that it is a length space and
either two points $x,y\in \wh{X}$ are connected by a geodesic in
$\wh{\mc{R}}$, or else for all sufficiently small $\eps$ a sequence of curves whose lengths limit to $d(x,y)$ eventually enters the $\eps$ neighborhood of a finite sequence of points $s_1,s_2,\dots,s_n\in\wh{\mc{S}}$. The geodesic between $x$ and $y$ then consists of the piecewise geodesic from $x$ to $s_1$ and $s_i$ to $s_{i+1}$ for $i=1,\dots,n$ and $s_n$ to $y$ since each segment lies in a single fundamental domain of $\mc{R}$ and is therefore globally length minimizing. Hence $\wh{X}$ is geodesic. Since each
component of $\wh{\mc{R}}$ has Riemannian sectional curvatures bounded above by
$-1$, it is enough to check the CAT(-1) condition in a neighborhood $U$ of a
completion point $p\in \wh{\mc{S}}$ not containing other points of
$\wh{\mc{S}}$. Since geodesics triangles straddling the two components of
$U\setminus p$ must pass through $p$, we need only check the conditions for
triangles with one vertex $p$ and the other two vertices in the same component
of $U\setminus p$. However, the Riemannian curvatures imply the CAT(-1)
condition for such triangles by a standard application of the Rauch comparison
theorem which still applies to the triangle away from $p$. This establishes the
claim.

The induced map $\wh{X}\to X$ is a ramified cover over $X$ with ramification locus
$\mc{S}\subset X$. Note that $\wh{X}$ is a nonproper metric space since small geodesic
spheres around points of $\wh{\mc{S}}$ will be homeomorphic to a disjoint union of the
(noncompact) universal cover of the closed cuspidal nilmanifolds on each side.

The graph of groups structure of $\pi_1(M)$ implies that it acts on $\wh{X}$ as well, but
edge stabilizer subgroups will fix corresponding points in $\wh{\mc{S}}$. Moreover, the
map $\wh{f}$ is equivariant with respect to this action on either side and the map $f$ is 
the equivariant quotient of $\wh{f}$ by this $\pi_1(M)$ action.

We now describe the family of natural maps $F_s:M\to X$ homotopic to $f$ and estimate the
Jacobian on $F_s^{-1}(\mc{R})$.

Let $U=f^{-1}(\mc{R})$ and $\til{U}\subset \til{M}$ be the lift of $U$ in the
universal cover of $M$. Let $g$ be any Riemannian metric on $M$ and for
all $s>h(g)$  consider the measure $\mu_x^s$ supported on $\til{U}$ absolutely
continuous with Lebesgue measure and with Radon-Nikodym derivative
$$\frac{d\mu_x^s}{d\vol_{\til{M}}}(z)=\chi_{\til{U}} e^{-s d(x,z)}.$$
Where $d$ is the distance on $\til{M}$. Note that the measure has finite total mass by
the condition on $s$.

Set $\sigma_x^s=\wh{f}_*\mu_x^s.$ Now consider the function $\mc{B}_{s,x}:\wh{X}\to \R$
given by
$$\mc{B}_{s,x}(y)=\int_{\wh{X}} d_{\wh{X}}(y,z)-d_{\wh{X}}(x,z)
d\sigma_x^s(z).$$

Note that the distance function on the CAT($-1$) space $\wh{\mc{R}}$ is smooth, and any
geodesic joining distinct components of $\wh{\mc{R}}$ passes through points of
$\wh{\mc{S}}$. It follows that $d_{\wh{X}}$ is smooth on $\wh{\mc{R}}\times \wh{\mc{R}}$.
In particular, since the support of $\sigma_x^s$ is $\wh{\mc{R}}$, $\mc{B}_{s,x}$ is
smooth on $\wh{\mc{R}}$. Moreover, it is strictly convex and we denote its unique minimum
point by $\op{Bar}(\sigma_x^s)$. We now set $\wh{F_s}(x)=\op{Bar}(\sigma_x^s)$. Since this
map is equivariant under the (isometric) actions of $\pi_1(M)$ on $\til{M}$ and $\wh{X}$,
it descends to a continuous map $F_s:M\to X$.

It is easy to verify that setting
$\wh{\Psi}_t(x)=\op{Bar}(t\delta_{\wh{f}(x)}+(1-t)\sigma_x^s)$ produces an explicit
equivariant homotopy from $\wh{F_s}=\wh{\Psi}_0$ to $\wh{f}=\wh{\Psi}_1$.

We now rely on a large scale local version of a key global estimate originally presented 
in \cite{Besson-Courtois-Gallot:95}.

\begin{Proposition}\label{prop:hyp_ball}
For all $s>h(g)$, the natural map $F_s:M\to X$ is $C^1$ on $F_s^{-1}(\mc{R})$. Moreover,
for any $\eps>0$ there is an $R_\eps>0$, depending only on $\eps,s,f$ and $M$, with
$\lim_{\eps\to 0}R_\eps=\infty$, such that for any component $V_i\subset\mc{R}$ and any
$x\in F_s^{-1}(V_i)$ whose ball $B(F_s(x),R_\eps)$ in $V_i$ is isometric to that of a rank
one locally symmetric space of with entropy $h_i$, then
$$\abs{\op{Jac} F_s(x)}\leq (1+\eps)\left(\frac{s}{h_i}\right)^n.$$
\end{Proposition}

A proof of the above proposition in the case when the comparison locally symmetric space 
is the three dimensional hyperbolic space, $\HH^3$, can be found in the proof of 
Proposition 4' in \cite{S}. 
In Appendix \ref{app:a}, we provide a different proof of proposition \ref{prop:hyp_ball} 
that works for all $n>2$ and when the pieces are modeled on any of the rank one symmetric 
spaces.

Now we construct a sequence of incomplete metrics $g_i^j$ on $V_i$ which are locally
symmetric on a set of diameter at least $j$, and whose remaining set has volume $\eps$.
We achieve this by taking a sequence of successively larger diameter locally symmetric
metrics on $M_j$ whose injectivity radius of each boundary component tends to $0$. For
each such metric we then close each end off as before with a negatively curved warped
product metric cusp of finite diameter as described above. Since the boundary components of 
the locally
symmetric portion may be chosen with arbitrarily small injectivity radius, both
the $n-1$ dimensional volume of the boundary component and the total volume of the 
non-locally symmetric finite cusps may be chosen to have total volume less than $\eps$.

Therefore, for each $\eps>0$ and constant $R_\eps$ from Proposition \ref{prop:hyp_ball} 
we choose $j>R_\eps$ and thus we can apply Proposition \ref{prop:hyp_ball} to the 
component $V_i$
with metric $g_i^j$. Putting these together we have metrics $g_j$ on all of $X$, and a
locally symmetric submanifold $(X_j,g_j)\subset (X,g_j)$ with $\vol(X-X_j,g_j)< k\eps$.
Hence,
$$\vol(X_j,g_j)\leq \int_{F_s^{-1}(X_j)}\abs{\op{Jac} F_s}d\vol_g\leq
(1+\eps)\sum_{i=1}^k \left(\frac{s}{h_i}\right)^n\vol(M_i,g),$$ and adding the remainder
yields:
$$\vol(X,g_j)\leq 2k\eps+\vol(X_j,g_j)\leq
2k\eps+(1+\eps)\sum_{i=1}^k \left(\frac{s}{h_i}\right)^n\vol(M_i,g).$$ Taking $\eps\to
0$, and thus $j\to\infty$, we obtain that $\vol(X_j,g_j)$ tends to
$\vol(M_{sym},g_{sym})$. Moreover, as will be shown in the proof of Theorem 
\ref{thm:minent} below, all of the pure pieces of a graph manifold must
be of the same rank one locally symmetric type so the $h_i's$ are all equal to a single
value $h_0$. Therefore, we obtain
$$\vol(M_{sym},g_{sym})\leq \frac{s^n}{h_0^n}\sum_{i=1}^k 
\vol(M_i,g)\leq\frac{s^n}{h_0^n}\vol(M,g).$$
Taking $s\to h(g)$ from above yields the result:
\[
h_0^n\vol(M_{sym},g_{sym})\leq h(g)^n\vol(M,g)
\]
\end{proof}

\begin{proof}[Proof of Theorem \ref{thm:minent}]

We first explain why, under the assumptions, all of the pure pieces must be of
the same type, that is, covered by the same locally symmetric space. For each
pure piece $Z_i=M_i$ of dimension $n$, each component of $\pa M_i$ is a compact
quotient of the same nilpotent group, either $\R^{n-1}$ in the real hyperbolic
case, or else one of three possible distinct 2-step groups of generalized Heisenberg
type: $N_1,N_2,N_3$ of dimensions $n-1=2k-1,4k-1,15$ respectively for some
$k\geq 2$. These quotients, with the exception of tori in the real hyperbolic
case, are irreducible and in particular not diffeomorphic to any nontrivial
product manifold. Neither can a quotient of $N_i$ be glued by an affine
diffeomorphism to a quotient of $N_j$ for $i\neq j$ since if their dimension are
the same, then the dimension of the centers differ. Therefore, since we assume
all of the pure pieces are locally symmetric, as soon as one of the pure pieces
is not real hyperbolic, then it can only be glued to a pure piece whose boundary
inherits its intrinsic metric from the same nilpotent group. In particular, this
implies that all the ends are covered by the same symmetric space. (If all the
pure pieces are real hyperbolic, then these may be glued to non pure pieces.)
Note that it is essential that we assume that all pure pieces are locally
symmetric. This is because a single pure piece that is not locally symmetric may
have distinct boundary components which are quotients of different nilpotent
groups.

The conclusion of Proposition \ref{prop:entropy} then gives the lower bound of the first
inequality of the entropy estimate. If we assume that the curvatures of $(M,g)$ are
scaled to have lower bound $-4$, then $h(g)\leq 2(n-1)$ by the Bishop-Gromov
comparison theorem. In summary,
$$h_0\Vol(M_{sym},g_{sym})^{\frac1n}\leq h(g)\vol(M,g)^{\frac 1n}\leq 2(n-1) 
\vol(M,g)^{\frac1n}.$$
Now observe that the quantity $h(g)\vol(M,g)^{\frac1n}$ is homothety invariant.
Hence, restricting to metrics within the given smooth class, we may express $h(M)$ as,
$$h(M)=\inf_g\set{h(g)\vol(M,g)^{\frac1n}}=\inf_g\set{h(g)\vol(M,g)^{\frac1n}\,\vert\, 
-4\leq K_g}.$$
We now produce a sequence of metrics $g_j$ on $M$ with $K_{g_j}\geq -4$ whose volumes
converge to $\vol(M_{sym},g_{sym})$. If entropy were Gromov-Hausdorff continuous then
this would show that $$h(M)=h_o\vol(M_{sym},g_{sym})^{\frac1n}.$$ However, since this is
unknown, we will be content using the bound $h(g_j)\leq 2(n-1)$.

To begin with, we note that boundary horospherical nilmanifolds of the locally
symmetric pieces are either flat in the real hyperbolic case, or else have a
generalized Heisenberg type with intrinsic curvatures lying in the interval
$[-3,1]$ when the locally symmetric metric has curvature bounds in $[-4,-1]$
(see \cite{N}). We now show that there is a sequence of metrics $g_j$ of
curvature bounded between $-4$ and $4$ on $M$ such that $\vol(M,g_j)\to
\vol(M_{sym},g_{sym})$.

Given a  $j\in \set{1,2,\dots}$, we construct the metrics $g_j$ on $M$ in two
stages as follows. Start with a locally symmetric metric on the interior of each
piece $M_i$ with curvatures in $[-4,-1]$. On each cusp, let $t$ be the
corresponding horospherical parameter (i.e. the negative of the corresponding
Busemann function), normalized so that quotient of the horosphere corresponding
to $t=0$ has intrinsic diameter $e^{-j-10}$. Let $\eps=\frac{1}{j}$ and $T=3j$,
and consider the portion of the cusp for $t\in [0,T]$. We will now modify the
locally symmetric metric on this portion of the cusp. (The compact manifold
$M_i$ is identified with the complement of the portions of all the cusps
corresponding to $t>0$.) Writing the cusp portion as $[0,T]\times Q$ for a
convex horospherical nilmanifold $Q$, the locally symmetric metric is a
generalized warped product metric $dt^2+f_0(z,t)dz^2$ where $dz^2$ represents a
metric on $Q$ which lifts to a left-invariant metric on the generalized
Heisenberg group $\tilde{Q}$ (or Euclidean space) and in the appropriate (left
invariant) coordinates on each tangent space of $Q$, $f_0(z,t)$ is diagonal with
entries $e^{-t}$ and one entry of $e^{-2t}$ corresponding to the intrinsic
dilatation. We replace this with a metric of the form $dt^2+f(z,t)dz^2$ where
$f(z,t)$ is a diagonal matrix with $e^{a(t)}$ and one entry of $e^{2a(t)}$, and
$a(t)$ will be constructed as follows.

Choose a smooth nonnegative function $[0,T]$ which satisfies $0< q(t)<\eps/2$ on $(0,T)$, vanishes to all orders at $t=0$ and $t=T$ and such that $\int_0^T q(t)dt=1.$ Such a function can be chosen since $T>\frac{2}{\eps}$. Since for all $t\in [0,T]$, $0\leq \int_0^t q(t)dt\leq 1$,  we have for either choice of $c\in \set{1,2}$,
\[
q(t)+c\left(-1+\int_0^t q(t)dt\right)< c+\eps/2.
\]
Now set $a''(t)=q(t)$ and suppose $a'(0)=-1$ and $a(0)=0$
, then if $h(t)=e^{c a(t)}$ for $c\in \set{1,2}$,
\[
0\leq\frac{h''(t)}{h(t)}=c(a''(t)+ca'(t))\leq c^2+c\eps/2\leq c^2+\eps.
\]
Similarly $0\geq \frac{h'(t)}{h(t)}=ca'(t)\geq -c$ for $t\in [0,T]$. The first and third equations of \eqref{eq:warped}, together with the comparison to the locally symmetric case where $a'(t)\equiv -1$, imply the sectional curvatures lie in $[-4-\eps,1]$.
Scaling the metric by the factor $\sqrt{1+\eps/4}$ yields a metric with curvatures in the interval $[-4,1]$ whose volume changes by a factor of $(1+\eps/4)^n$. The resulting metric transitions from the nilpotent metric on the slice $\set{0}\times Q$ inherited from the locally symmetric one to one on $\set{1}\times Q$ that is extrinsically totally geodesic. 

The above construction subsumes the case of real hyperbolic pieces as well. In that case the boundary components in the induced metric above are totally geodesic flat manifolds, and there is no $c=2$ case. Hence the resulting curvature of these metrics lie in the interval $[-1-\eps,0]$. 

The non-pure pieces combine together to form connected components which will
glue to the symmetric pieces. Each of these components admits a sequence of fibered metrics, indexed by $j\in \N$, with structure of a nilmanifold bundle over a negatively curved manifold and such that the curvatures lie in $[-1,1]$ and the volume is at most $e^{-jn}$. In this sequence the entire center of the nil-fibers is uniformly collapsed. Without loss of generality we may
assume that the boundary manifolds of these components remain totally geodesic
nilmanifolds covered by Euclidean space or else generalized Heisenberg type
groups and that as in Lemma \ref{lemma:bilip}, the attaching diffeomorphism on the boundary components are uniformly Bilipschitz independent of $j$. 

The gluing diffeomorphisms of the pure pieces to other pure or non-pure pieces
will be isotopic to affine maps which are induced from an automorphism of the
underlying nilpotent Lie group. Hence, even if the boundaries in these metrics
were totally geodesic, we still cannot directly join these metrics on the pieces
together unless the affine map is an isometry. Instead, we will proceed
similarly to the proof of Theorem \ref{thm:zero}.

If $N=\partial Z_1=\partial Z_2$ is the boundary nilmanifold with nilpotent
universal cover $\widetilde{N}$ between two paired boundary components $Z_1$ and
$Z_2$, then we will build a patch metric on $A=[-R,R]\times N$, for some fixed
$R>4$, so that the manifolds may be glued isometrically to the boundaries of the
patch in such a way that the resulting manifold will be diffeomorphic to the
original higher graph manifold.

On $A$ we define the metric to be $f(z,t)=\alpha(t)f_1(z)$ for $t\in [-R,-1]$ and $f(z,t)=\alpha(-t)f_1(z)$ for $t\in [1,R]$ where $f_i$ is the boundary metric of the boundary component of $Z_i$, and $\alpha:[-R,-1]\to \R$ is smooth with $\alpha(-R)=1$ and $\alpha'(t)\in [0,1]$ vanishing to all orders at $t=-R$ and $t=-1$, $-2\leq \alpha''(t)\leq 2$ for $t\in [-R,-1]$ and $\alpha'(t)\equiv -1$ for $t\in [-R+2,-3]$. Since $-2\leq \frac{\alpha''(t)}{\alpha(t)}\leq 2$ and $0\leq \frac{\alpha'(t)}{\alpha(t)}\leq 1$, from \eqref{eq:warped} we see that the curvatures of this section of $A$ lie in $[-4,4]$. However, we note that scaling factor is no longer the Lie algebra dilatation, but a straight scaling of the fiber metric. In particular the metric on the copies of $Q$ at $t=1$ and $t=-1$ have curvatures lying in $[-\frac{3}{R-2},\frac{1}{R-2}]$. Their diameters are the same and so from Lemma \ref{lemma:bilip} the affine diffeomorphism between them is $C$-bilipschitz and the metrics are $C$-relatively bounded up to first order for some $C>0$. On the remaining part $[-1,1]\times Q$ we put the metric $dt^2+(1-t)/2\rho_1+(t+1)/2 \rho_2$ where $\rho_1$ and $\rho_2$ are the metrics at $t=-1$ and $t=1$ respectively. 
From Lemma \ref{lem:curv_of_sums} we see that the curvatures are bounded in $[-4,4]$ for sufficiently large, but independent of $j$, choice of $R$. (Note in the proof of Lemma \ref{lem:curv_of_sums} that the remainder $C_{ijkl}$ tensor components also scale with the metric.)
    
%and constant in $t$ for $t\in [2,3]$ on $Z_2$.
%We will now build collapsing metrics on $A$ whose collapsing central 
%directions in each fiber coincide with the collapsing directions in the portion of the 
%fibered metric $f(z,t)=f(z)$ where they are averaged together using bump functions for 
%$t\in [-3,-2]$ and $t\in [2,3]$. The existence of such metrics follows a similar procedure to the proof of Theorem \ref{thm:zero} above (and ultimately Fukaya's ideas). We assume the fiber metrics here are sufficiently small 
%so that the volume of all of $A$ is less than $e^{-j}$.

From crude estimates from the lengths and diameter bounds given above, and accounting for the rescaling, the volume of $g_j$ on each cusp, away from the portion of $(M,g_j)$ which carries a locally symmetric
metric, is bounded above by $(1+\frac{1}{4j})^n(4j+2(R)^{n+1})e^{-jn}$. Since the connected components of attached nonpure pieces have volumes at most $e^{-jn}$, taking a limit as $j\to\infty$, the
sequence of metrics $(M,g_j)$ has volumes converging to $\Vol(M_{sym},g_{sym})$. Since the curvatures lie in $[-4,4]$, the entropy estimate mentioned at
the beginning gives us the upper bound of the first line of inequalities 
\eqref{Thm3-ineq1}.

For the second statement of the theorem, namely the inequalities
\eqref{Thm3-ineq2}, if $(M,g)$ has curvature bounds $K_g\geq -1$, then $h_0\geq
n-1\geq h(g)$ and from the lower bound of the first statement we have,
$$\frac{h_0}{n-1}\vol(M_{sym},g_{sym})^{\frac1n}\leq \frac{h_0}{h(g)} 
\vol(M_{sym},g_{sym})^{\frac1n}\leq \frac{h(M)}{h(g)}\leq \vol(M,g)^{\frac1n}.$$
The lowermost bound then follows by taking infimums in $g$. (Incidentally, this
inequality will be strict unless the pure pieces are real hyperbolic.) For the
uppermost bound, we note that the metrics $g_j$ constructed above, when scaled
to have curvatures bounds $-1\leq K_g\leq 1$, have volumes tending to $2^n
\vol(M_{sym},g_{sym})$. Finally, we observe that $\Vol_{K}(M)\leq \Minvol(M)$
follows in general from their definitions.
\end{proof}

\begin{Remark}
In \cite{FLS} the authors construct some simple and not-so simple examples of
higher graph manifolds, though without pure pieces, which admit no CAT$(0)$ metric. They 
do, however, show in their theorem 2.8 that graph manifolds with only hyperbolic pieces 
(purely hyperbolic graph manifolds) indeed support a nonpositively curved metric. 
\end{Remark}

\begin{proof}[Proof of Theorem \ref{thm:minvol}]
Observe that under the hypothesis of hyperbolic pure pieces, the metrics $g_j$ on $M$
constructed in the proof of Theorem \ref{thm:minent} have curvatures lying in the
interval $[-1,0]$, since the boundary {infranilmanifolds} in this case are just flat 
manifolds, though not necessarily tori.

Notice that $\vol(M,g_j)$ tends to $\vol(M_{hyp},g_{hyp})$. Combining this with the
inequalities in equation \ref{Thm3-ineq2} of Theorem \ref{thm:minent} completes the proof.
\end{proof}

Before we proceed with the proofs of the rest of our results, let us include a proof here 
of the fact that higher graph  manifolds are aspherical. To this end we also include the 
next definition, following \cite{FLS}:

\begin{Definition}
The boundaries of the pieces $Z_{i}$ that are identified together in a higher graph 
manifold $M$ will be called the {\em internal walls} of $M$.
\end{Definition}

\begin{Lemma}\label{asphericity} 
If $M$ is a higher graph manifold (possibly with boundary), then $M$ is aspherical.
\end{Lemma}

\begin{proof}
We will work by induction on the number of internal walls $c$ of $M$. If $c=0$ then $M=Z$ 
for some infranil bundle $Z$ over a closed, negatively curved base. It follows from the 
homotopy exact sequence for the bundle $Z$ that $M$ is aspherical in this case, 
establishing the base case for our inductive argument.

Assume $c>0$, and that the result holds for higher graph manifolds with strictly less 
than $c$ internal walls. Cut open $M$ along an arbitrary internal wall $W$. Our inductive 
hypothesis implies that now $M$ is obtained by gluing one or two (depending on whether 
$W$ separates $M$ or not) aspherical spaces. Since the inclusion of $W$ in the piece(s) 
in $M$ it belongs to is $\pi_{1}$--injective, it follows from a classical result of 
Whitehead \cite{W39} that $M$ is aspherical. \end{proof}

We thank J.-F. Lafont for suggesting the strategy of the following proof and pointing us 
to the relevant references. The concept of {\it asymptotic dimension} of a group was 
introduced by M. Gromov in \cite{Gro2}:

\begin{Definition}
 Let $X$ be a metric space. We say that the asymptotic dimension of $X$ does not exceed 
 $n$ and write ${\rm asdim} X \leq n$ provided that for every uniformly bounded open 
 cover $V$ of $X$ there is a uniformly bounded open cover $U$ of $X$ of multiplicity 
 $\leq n+1$ so that $V$ refines $U$. We write ${\rm asdim} X = n$ if it is true that 
 ${\rm asdim} X \leq n $ and ${\rm asdim} X\nleq n-1$.
\end{Definition}

We refer the reader to \cite{BD} for examples of spaces and groups with finite asymptotic 
dimension.

\begin{proof}[Proof of Theorem \ref{thm:psc}]

Let $Z_{i}$ be a piece of the higher graph manifold $M$. To begin with we claim that
$\pi_1(Z_{i})$ has finite asymptotic dimension.

The infranilpotent factor $N_{i}$ is finitely covered by a nilmanifold $N'_{i}$. It
follows from results of Bell and Danishnikov, Corollaries 54 \& 68 in \cite{BD}, that the
asymptotic dimension of $\pi_1(N_{i})$ is finite.

The main result of Osin in \cite{O} shows that in this situation the asymptotic dimension
of $\pi_1(M_{i})$ is finite.
The total space $Z_{i}$ of the fibration $N_{i}\hookrightarrow Z_{i}\to M_{i}$ fits into 
a short exact sequence
$$ 1\to  \pi_1(N_{i})\to \pi_1(Z_{i})\to \pi_1(M_{i}) \to 1.  $$
A theorem that deals with the asymptotic dimension of these types of extensions of groups 
has been shown by Bell--Dranishnikov, see Theorem 63 of \cite{BD}.
They show that the asymptotic dimension of $\pi_1(Z_{i})$ is bounded from above by the 
sum of the asymptotic dimensions of $\pi_1(N_{i})$ and $\pi_1(M_{i})$. Therefore 
$\pi_1(Z_{i})$ has finite asymptotic dimension.

Here we will again use that $\pi_1(M)$ is  a graph of groups, with vertex groups
$\pi_1(Z_{i})$. Each vertex group $\pi_1(Z_{i})$ has finite asymptotic dimension, so
Theorem 77 in \cite{BD} implies that  $\pi_1(M)$ has finite asymptotic dimension, as
claimed.

In lemma \ref{asphericity} we have proven that
$M$ is aspherical. It is now a consequence of the results of G. Yu in \cite{Y}, that, 
since $M$ is
aspherical and $\pi_1(M)$ has finite asymptotic dimension, the manifold $M$ does not
admit any smooth metric of positive scalar curvature. \end{proof}

Next, we present the coarse Baum-Connes conjecture for these groups.
\begin{proof}[Proof of Proposition \ref{cor:bc}]
The main result of \cite{Y} states that a finitely generated group $G$ satisfies the 
coarse Baum-Connes conjecture if its classifying space $BG$ has the homotopy type of a 
finite CW-complex and if $G$ has finite asymptotic dimension with respect to any left 
invariant metric. This last condition we verified in Theorem \ref{thm:psc}. 

The asphericity of $M$ was shown in Lemma \ref{asphericity}, following the 
proof of the simpler case treated in \cite{FLS}. Hence the finite CW-complex  $M$ serves 
as a $B\pi_1(M)$.
\end{proof}

\begin{proof}[Proof of Corollary \ref{cor:yam}]
If a smooth compact manifold $M$ admits a volume collapsing sequence of metrics with
bounded curvature then $ \mathcal{Y}(M) \geq 0$  \cite{Sc}. On the other hand, a well
known fact about the Yamabe invariant is that $\mathcal{Y}(M) >0 $ if and only if $M$
admits a smooth metric of positive sectional curvature, (we refer again to \cite{Sc}).
Therefore Theorems \ref{thm:zero} and \ref{thm:psc}  imply $\mathcal{Y}(M) =0$.
\end{proof}

\appendix
\section{Proof of Proposition \ref{prop:hyp_ball}}\label{app:a}

For convenience we restate the Proposition here. 

{
\renewcommand{\theTheorem}{\ref{prop:hyp_ball}}
\begin{Proposition}
For all $s>h(g)$, the natural map $F_s:M\to X$ is $C^1$ on $F_s^{-1}(\mc{R})$. Moreover,
for any $\eps>0$ there is an $R_\eps>0$, depending only on $\eps,s,f$ and $M$, with 
$\lim_{\eps\to 0}R_\eps=\infty$, such that
for any component $V_i\subset\mc{R}$ and any $x\in F_s^{-1}(V_i)$ whose ball
$B(F_s(x),R_\eps)$ in $V_i$ is isometric to that of a rank one locally symmetric space of 
with entropy $h_i$, then
$$\abs{\op{Jac} F_s(x)}\leq (1+\eps)\left(\frac{s}{h_i}\right)^n.$$
\end{Proposition}
\addtocounter{Theorem}{-1}
}
This proposition and its proof are motivated by Proposition 4' in \cite{S}. However, we 
will use a somewhat different approach to handle the more general setting. We will 
assume some familiarity with the prototype proof along these lines found in 
\cite{Besson-Courtois-Gallot:95}.

\begin{proof}
Let $\wh{\mc{R}}\subset \wh{X}$ denote the lift of $\mc{R}\subset X$, which is the set of 
regular manifold points in $\wh{X}$.  The pre-image ${F}_{s}^{-1}{(\mc{R})}$ of the regular 
set is an open subset of $M$. 

We will first show that ${F}_{{s}}$ is $C^{1}$ when restricted to 
${F}_{s}^{-1}{(\mc{R})}$.

Since $\wh{X}$ is a CAT(-1) space, it is uniquely geodesic and geodesics vary continuously in 
the endpoints. However, such a metric in general may not be smooth even when restricted 
to the portion of the space admitting a smooth structure compatible with the topology 
induced by the metric. Nevertheless, in our case the singular set consists of isolated 
points, and all geodesics from any fixed $p\in \wh{X}$ to any point in one component of the 
regular set have a common fixed segment containing all of the singular points that the 
geodesic passes through. Hence $d_{\wh{X}}(p,\cdot)$ is smooth on the regular set for any 
$p\in \wh{X}$.

Since ${\sigma}_{x}^{s}$  is supported on $\wh{\mc{R}}$, for $x\in \til{M}$ the 
function $\mc{B}_{s,x}:\wh{X}\to \R$ is smooth on $\wh{\mc{R}}$ with gradient

\[
\nabla_y{\mc{B}}_{s,x}=\int_{\wh{X}}\nabla_y \phi\, d{\sigma}_{x}^{s}(z),
\]
where $\phi$ is the 1-Lipschitz function $\phi(y)=d_{\wh{X}}(y, z)$ whose gradient 
should be understood in the weak sense. In particular, for 
$x\in \wh{F_{s}}^{-1}(\wh{\mc{R}})$, we have the defining equation
\[
\nabla_{{\wh{F}}_{s}(x)}{\mc{B}}_{s,x}=0.
\]

Setting $g(x)=\sum_{q\in \wh{f}^{-1}(z)} d_{\til{M}}(x,q)$ we may covariantly 
differentiate the vanishing field $\nabla_{{\wh{F}}_{s}(x)}{\mc{B}}_{s,x}$ in the 
direction $u\in T_x\til{M}$ yielding, 

\begin{align*}
D_u \nabla_{\wh{F_s}(x)} {\mc{B}}_{s,x}
&=\int_{\wh{X}} (D_\cdot\nabla_{y}\phi)|_{y=\wh{F_s}(x)}  \of 
d_x\wh{F_s}(u)d{\sigma}_{x}^{s}(z)\\
&\hspace{2cm} + \int_{\wh{X}}\nabla_{\wh{F_s}(x)} \phi 
\tensor d_x \left( \frac 
{d{\sigma}_{x}^{s}}{d{\sigma}_{p}^{s}}(z)\right)(u)\, 
d{\sigma}_{p}^{s}(z)\\
&= \left(\int_{\wh{X}} 
(D_\cdot\nabla_{y}\phi)|_{y=\wh{F_s}(x)}d{\sigma}_{x}^{s}(z)\right)\of 
d_x\wh{F_s}(u)\\
&\hspace{2cm} - s\int_{\wh{X}}\nabla_{\wh{F_s}(x)} \phi \tensor d_x g(u)\, 
d{\sigma}_{x}^{s}(z).
\end{align*}

Here we have assumed that the associated objects exist at least weakly in $L^1$. We will 
address this shortly.

Setting
\[
K_x^s=\int_{\wh{X}} 
(D_\cdot\nabla_{y}\phi)|_{y=\wh{F_s}(x)}d{\sigma}_{x}^{s}(z) \quad\text{and}\quad
A_x^s(u)=\int_{\wh{X}}\nabla_{\wh{F_s}(x)} \phi \tensor d_x g(u)\, 
d{\sigma}_{x}^{s}(z),
\]
we can use that $0=D_u \nabla_{\wh{F_s}(x)} {\mc{B}}_{s,x}$ to solve for $d_x\wh{F_s}(u)$ 
giving
\[
d_x\wh{F_s}=s (K_x^s)^{-1}\of A_x^s.
\]

Here the $(1,1)$-form $D\nabla_y \phi=II_y\oplus 0$ where $II_y$ is the second fundamental 
form at the point $y$ of the sphere of radius  $\phi(y)$ centered at $z$ operating on its 
tangent space, which is then extended to be 0 in the normal $\nabla_y \phi$ direction. 
Since the spheres in any CAT($-1$) space are strictly convex, and in our case smooth, the form
$D\nabla_y \phi$ is positive definite except in the null direction $\nabla_y \phi$.

Since $\sigma_x^s$ is absolutely continuous with respecto to the Lebesgue measure and the singular set 
in our case is discrete, we will show that the integrand of $K_x^s$ is almost everywhere 
bounded away from $0$ from below and makes sense. However, a priori it could be singular 
due to the contributions near the singular set. Nevertheless, because of the negative 
Riemannian curvature the integrand is well defined for every $z\in \wh{X}$ and $y\in 
\wh{\mc{R}}$. (Recall the singular set has measure zero and consists of a discrete set of 
points.)

To understand the integrand, we observe from that in constant cuvature $-k^2\leq -1$ the
solutions to the Ricatti equation imply that the second fundamental form at any point $y$
of a sphere of radius $t$ is $II_y=k\coth(k t)\geq \coth(t)$. By the standard comparison
estimate (see e.g. \cite{CE}), for the metric on the regular part of $\wh{X}$ coming from
the warped product construction on the symmetric space, even when $z$ is a singular point,
we always have $II_y \geq \coth(d(y,z))$ since its sectional curvatures are strictly less
than $-1$. Even when $z$ is a singular point, for $y$ restricted to a region where the
curvature is always greater than $-k^2$, for instance the complement of some
$\eps(k)$-neighborhood of the singular points, then we also have the reverse inequality
$k\coth(k t) \geq II_y$.

Nevertheless, the operator $K_x^s:T_{\wh{F_s}(x)}\wh{X}\to T_{\wh{F_s}(x)}\wh{X}$ may 
have an infinite eigenvalue if an eigenvalue of $II_{y}$ grows sufficiently rapidly 
as $y$ approaches $z$. This would lead to a nonzero element in the kernel of $d_x\wh{F_s}$.

Since $\sigma_x^s$ is non atomic, the average of the positive-semidefinite forms
$D\nabla_y \phi$ each admitting a single null eigen-direction, which is distinct for each
$y$, will be strictly positive definite. Hence $K_x^s$ has full rank at each $x$ in the
stated domain. In particular, the implicit function theorem implies that $\wh{F_s}$ is
$C^1$ on ${F}_{s}^{-1}{(\mc{R})}$.

By Cauchy-Schwarz applied to bilinear forms, we may write,
\[
(A_x^s)^2\leq H_x^s \of B_x^s,
\]
where
\[
H_x^s=\int_{\wh{X}}\nabla_z \phi \tensor d_zf \, d{\sigma}_{x}^{s}(z) \quad\text{and}\quad 
B_x^s=\int_{\wh{X}}\nabla_z g\tensor d_z g\, d{\sigma}_{x}^{s}(z).
\]
The determinant of $B_x^s$ can be estimated by noting that the trace of the integrand is 
one, and that among such matrices, the determinant will be maximal when 
$B_x^s=\frac{1}{n}I$. Therefore $\det B_x^s\leq \frac{1}{(n_i)^{n_i}}$. Consequently,

\[
\left(\Jac \wh{F_s}(x)\right)^2\leq \left(\frac{s^{2}}{n_i}\right)^{n_i} \frac{\det 
H_x^s}{(\det K_x^s)^2}.
\]

Now we estimate $K_x^s$ in terms of the corresponding object in the symmetric space.
Recall that at each point of  $\mc{R}$, the metric on $X$ and that of the symmetric space
$\HK$ agree on some ball. Assume for the moment that $B(p,R)$ is contained in this ball 
and we consider any point $z\in S(p,R)$. If $\bar{\phi}(\bar z)=d_{\HK}(\bar{p},\bar z)$ 
represents the corresponding distance function in the symmetric space $\HK$ to a comparison 
point $\bar{p}\in \HK$ with $d(\bar{p},\bar z)=d(p,z)$, then
\[
D\nabla_z \phi(z) =D\nabla_{\bar z} \bar{\phi}(\bar z).
\]

The expression on the right hand side is
\[
D\nabla_z\bar{\phi}(\bar z)_{|_{(\grad \bar{\phi})^\perp}}=\sqrt{-R_{\bar 
z}}\coth(\bar{\phi}(\bar 
z)\sqrt{-R_{\bar z}}).
\]
Here $R_{\bar{z}}$ is the (1,1)-form dual to the curavature tensor $R(\grad 
\bar{\phi},\cdot,\grad,\bar{\phi},\cdot)$ at the point $\bar z$ twice contracted in the 
direction of $\grad \bar{\phi}$. This formula 
follows by explicitly solving the Ricatti equation using the fact 
that the curvature tensor is parallel in $\HK$. 

If the base field $\mb{K}$ has real dimension $d$, then we have $d-1$ complex structures 
$J_i:T\HK\to T\HK$, each satisfying $J_i^{-1}=-J_i$ which allows us to explicitly write,
\[
\sqrt{-R_{\bar z}}=I-E_{\bar z}-\sum_{i=1}^{d-1} J_i E_{\bar z} {J_i}_{|_{(\grad 
\bar{\phi})^\perp}}
\]
where $E_{\bar z}$ is the $(1,1)$-form $\grad_{\bar z} \bar{\phi}\tensor d_{\bar 
z}\bar{\phi}$. 
In the direction of $\grad \bar{\phi}$, we have that $D\nabla_{\bar z}\bar{\phi}(\bar z)=0$. 
So we 
may write,
\[
D\nabla_{\bar z}\bar{\phi}=(I-E_{\bar z}-\sum_{i=1}^{d-1} J_i E_{\bar z} 
{J_i})\coth\left(\bar{\phi}(\bar z)(I-E_{\bar z}-\sum_{i=1}^{d-1} J_i E_{\bar z} {J_i})\right).
\]
On the other hand, since $\coth(t)\geq 1$ for $t>0$ we have
\[
\coth\left(\bar{\phi}(\bar z)(I-E_{\bar z}-\sum_{i=1}^{d-1} J_i E_{\bar z} {J_i})\right)\geq I
\]
as positive definite symmetric two forms. For any number $R>0$ for which 
$B(\wh{F_s}(x),R)\subset \wh{X}$ is isometric to the ball $B(\bar{p},R)\subset \HK$ we can 
also construct a (strictly smaller) comparison measure $\bar{\sigma}_{x}^s$ on $\HK$ by 
setting it to be $\sigma_x^s$ on the ball $B(\bar{p},R)$ and $0$ outside of this ball. 

The maps $J_i$ commute with the action of the maximal compact subgroup $K<{\rm 
Isom}(\HK)$. Hence after integrating we have,
\begin{align*}
\bar{K}_{x}^s &:= \int_{\HK} D\nabla_z\bar{\phi}(\bar z) d\bar{\sigma}_{x}^s(\bar z)\geq 
\int_{\HK} 
I-E_{\bar z}-\sum_{i=1}^{d-1} J_i E_{\bar z}J_i \ d\bar{\sigma}_{x}^s(\bar z)\\ 
&= I-\bar{H}_{x}^s-\sum_{i=1}^{d-1} J_i \bar{H}_{x}^s J_i,
\end{align*}
where we have defined
\[
\bar{H}_{x}^s:= \int_{\HK} E_{\bar z} d\bar{\sigma}_{x}^s(\bar z).
\]

Since $K_x^s\geq\bar{K}_x^s$, we can combine these results to obtain,
\[
\left(\Jac \wh{F_s}(x)\right)^2\leq \left(\frac{s^{2}}{n_i}\right)^{n_i} \frac{\det 
H_x^s}{\det(I-\bar{H}_{x}^s-\sum_{i=1}^{d-1} J_i \bar{H}_{x}^s J_i)^2}.
\]
Now the mass of the measure $\norm{\bar{\sigma}_x^s}$ tends to $\norm{\sigma_x^s}$ as 
$R\to\infty$.  Since the measure $\sigma_x^s$ is fully supported on 
$B(\bar{p},R)$, $\bar{H}_x^s$ is (strictly) positive definite. 

Writing $H_x^s=\bar{H}_x^s+Q$ where 
\[
Q=\int_{X\setminus B(\wh{F_s}(x),R_\eps)} \nabla_z \phi \tensor d_zf d\sigma_{x}^s(z).
\] 
We note that since $\norm{\nabla_z \phi \tensor d_zf}$ is bounded by $1$, we have $Q$ tends 
to the $0$ form as $R\to \infty$. In particular for any $\eps>0$ there is a smallest 
$R_\eps$ such that $\det H_x^s \leq (1+\eps)\det \bar{H}_x^s$.

Hence 
\[
\left(\Jac \wh{F_s}(x)\right)^2\leq (1+\eps)\left(\frac{s^{2}}{n_i}\right)^{n_i} 
\frac{\det \bar{H}_x^s}{\det(I-\bar{H}_{x}^s-\sum_{i=1}^{d-1} J_i \bar{H}_{x}^s J_i)^2}.
\]

Now we observe that $R_\eps$ also depends continuously on $\sigma_x^s$ which depends
continuously on $s$, $x$,  the metric on $M$ and the underlying map $f$. So $R_\eps$
depends on these quantities and the metric on $X$. However, for a fixed metric on $M$ and
map $f$, we observe that the set $\til{L}\subset \til{M}$ consisting of $x\in
\wh{F_s}^{-1}(\wh{\mc{R}})$ such that $\wh{F_s}(x)$ has an open ball of radius one
isometric to the open ball of radius one in the symmetric space $\HK$ has pre-compact
quotient $L\subset M$. So for $R>1$ the equivariance of $\wh{f}$ implies that the fraction
of the total mass of $\sigma_x^s$ lying in a ball of radius $R$ that is isometric to a
ball of radius $R$ in $\HK$ can be bounded uniformly from below for all $x\in \til{L}$.
Hence $R_\eps$ can be chosen independent of $x$ and the metric on $X$. Consequently
$R_\eps$ only depends on $\eps, f$ and the metric on $M$.
  
We observe that since the measure $\sigma_x^s$ is nonatomic, the 2-form $\bar{H}_x^s$ is 
strictly positive definite. The proposition is then completed by the following lemma. \end{proof}

\begin{Lemma}[Proposition B.1 and B.5 of \cite{Besson-Courtois-Gallot:95}]
For all $n\times n$ positive definite matrices $A$ with trace one, and orthogonal 
matrices $J_1,\dots, J_{d-1}$ with $J_i^2=-I$ we have
\[
\frac{\det A}{\det(I-A -\sum_{i=1}^{d-1} J_i A J_i)^2} \leq 
\left(\frac{n}{(n+d-2)^{2}}\right)^n 
\]
with equality if and only if $A=\frac{1}{n}I$.
\end{Lemma}

\end{document}